\documentstyle[11pt,amssymb,amstex]{amsart}

\numberwithin{equation}{section}

\newcommand{\bea}{\begin{eqnarray}}
\newcommand{\eea}{\end{eqnarray}}
\newcommand{\ba}{\begin{array}}
\newcommand{\ea}{\end{array}}
\newcommand{\edc}{\end{document}}
\newcommand{\bc}{\begin{center}}
\newcommand{\ec}{\end{center}}
\newcommand{\be}{\begin{equation}}
\newcommand{\ee}{\end{equation}}

\newcommand{\dsf}{\displaystyle\frac}

\def\bc{{\mathbb C}}

\def\bn{{\mathbb N}}

\def\bq{{\mathbb Q}}
\def\br{{\mathbb R}}

\def\bz{{\mathbb Z}}

\def\a{\alpha}
\def\b{\beta}
\def\g{\gamma}

\def\e{\epsilon}

\def\h{{\mathbf{h}}}
\def\z{{\mathbf{s}}}

\newtheorem{thm}{Theorem}[section]
\newtheorem{lem}[thm]{Lemma}

\begin{document}
\small

\title[On a recursive equation]
{On a recursive equation over $p$-adic field} \thanks{On leave
from Department of Mechanics and Mathematics, National University
of Uzbekistan, Tashkent, 700174, Uzbekistan}
\author{Farrukh Mukhamedov}
\address{Farrukh Mukhamedov\\
 Departamento de Fisica\\
Universidade de Aveiro\\ Campus Universit\'{a}rio de Santiago\\
3810-193 Aveiro, Portugal} \email{{\tt far75m@@yandex.ru} {\tt
farruh@@fis.ua.pt}}

\begin{abstract}
In the paper we completely describe the set of all solutions of a
recursive equation, arising from the  Bethe lattice models over
$p$-adic numbers. \vskip 0.3cm \noindent {\it
Mathematics Subject Classification}: 46S10, 12J12.\\
{\it Key words}: recursive equation, $p$-adic numbers, model,
fixed point, unique.
\end{abstract}

\maketitle

\section{introduction}

The $p$-adic numbers were introduced by K. Hensel. For about a
century after the inventing of the $p$-adic numbers, they were
mainly considered objects of pure mathematics. After discovering
that the physics of certain models could be based on the idea that
the structure of space-time for very short distances might
conveniently be described in terms of $p$-adic numbers, many
applications of such numbers in theoretical physics have been
proposed in works (see for example,
\cite{FW},\cite{Kh1},\cite{MP},\cite{VVZ}).  A number of $p$-adic
models in physics cannot be described using ordinary probability
theory based on the Kolmogorov axioms, therefore $p$-adic
probability models were investigated in \cite{Kh2}. Such models
appear to provide the probabilistic interpretation of $p$-adic
valued wave functions and string amplitudes in the framework of
$p$-adic theoretical physics (see \cite{Kh1,Kh2}). Using such a
$p$-adic probability approach in \cite{MR1,MR2} we have developed
a theory of statistical mechanics in the context of that theory of
probability. In statistical mechanics, due to their solvable
character \cite{B}, Bethe lattice models and models on random
tree-like graphs \cite{DM}, play a central role. In the
calculation of the average magnetization of these models naturally
appears the following recursive equation
\begin{equation}\label{eq}
h_n=\log\bigg[\left(\dsf{e^{\a_{n+1}+h_{n+1}}+e^{\b_{n+1}}}{e^{\g_{n+1}}+e^{h_{n+1}}}\right)
\left(\dsf{e^{\a_{n+2}+h_{n+2}}+e^{\b_{n+2}}}{e^{\g_{n+2}}+e^{h_{n+2}}}\right)\bigg],
\ \ n\in\bn,
\end{equation}
here $\a_k,\b_k,\g_k\in\br$ and $\{h_n\}$ is an unknown sequence
of real numbers, the solution of such a equation describes
corresponding Gibbs measures. On certain conditions depending on
parameters $\a_k,\b_k,\g_k$ the equation (\ref{eq}) has infinitely
many solutions \cite{BG},\cite{G} (i.e. over real numbers), but
all the set of solutions is still not described.

In the present paper we are going to completely describe the set
of solutions of (\ref{eq}) over $p$-adic numbers. Our hope is that
this allows us further understanding of the Bethe lattice models
and models on random tree-like graphs over $p$-adic numbers.
Moreover, we think that the obtained result will give certain
information on random rational $p$-adic dynamical systems
generated by linear-fractional functions. Note that some
applications of
$p$-adic dynamical systems to some biological and physical systems have been proposed in \cite{KhN}.\\

Throughout the paper $p$ will be a fixed prime number greater than
3, i.e. $p\geq 3$. Every rational number $x\neq 0$ can be
represented in the form $x=p^r\dsf{n}{m}$, where $r,n\in\bz$, $m$
is a positive integer, $(p,n)=1$, $(p,m)=1$. The $p$-adic norm of
$x$ is given by
$$
|x|_p=\left\{ \ba{ll}
p^{-r} & \ \textrm{ for $x\neq 0$}\\
0 &\ \textrm{ for $x=0$}.\\
\ea \right.
$$

It satisfies the following strong triangle inequality
$$
|x+y|_p\leq\max\{|x|_p,|y|_p\},
$$
this is a non-Archimedean norm.

The completion of the field of rational numbers $\bq$  with
respect to the $p$-adic norm is called  {\it $p$-adic field} and
it is denoted by $\bq_p$.

Given $a\in \bq_p$ and  $r>0$ put
$$
B(a,r)=\{x\in \bq_p : |x-a|_p< r\}, \ \ S(a,r)=\{x\in \bq_p :
|x-a|_p=r\}.
$$
The {\it $p$-adic logarithm} is defined by the series
$$
\log_p(x)=\log_p(1+(x-1))=\sum_{n=1}^{\infty}(-1)^{n+1}\dsf{(x-1)^n}{n},
$$
which converges for $x\in B(1,1)$; the {\it $p$-adic exponential}
is defined by
$$
\exp_p(x)=\sum_{n=0}^{\infty}\dsf{x^n}{n!},
$$
which converges for $x\in B(0,p^{-1/(p-1)})$.

\begin{lem}\label{1}\cite{Kl} Let $x\in B(0,p^{-1/(p-1)})$
then we have
\begin{eqnarray}\label{exp}
&& |\exp_p(x)|_p=1,\ \ \ |\exp_p(x)-1|_p=|x|_p,
\ \ \ |\log_p(1+x)|_p=|x|_p \\
\label{el} && \log_p(\exp_p(x))=x, \ \ \exp_p(\log_p(1+x))=1+x.
\end{eqnarray}
\end{lem}

\section{Main result}
In this section  we are going to completely describe the set of
solutions of (\ref{eq}) over $p$-adic numbers. Let us first
formulate the problem in a $p$-adic setting. Namely, we are
interested to find all solutions of the following recursive
equation
\begin{equation}\label{eq1}
h_n=\log_p\bigg[\left(\dsf{\exp_p(\a_{n+1}+h_{n+1})+\exp_p(\b_{n+1})}{\exp_p(\g_{n+1})+\exp_p(h_{n+1})}\right)
\left(\dsf{\exp_p(\a_{n+2}+h_{n+2})+\exp_p(\b_{n+2})}{\exp_p(\g_{n+2})+\exp_p(h_{n+2})}\right)\bigg],
\ \ n\in\bn,
\end{equation}
here $\a_k,\b_k,\g_k\in\bq_p$, such that
\begin{equation}\label{con}
|\a_k|_p<p^{-1/(p-1)}, \ \ |\b_k|_p<p^{-1/(p-1)}, \ \
|\g_k|_p<p^{-1/(p-1)},
\end{equation}
and $h_k\in\bq_p$ also should satisfy $|h_k|_p<p^{-1/(p-1)}$ for
every $k\in\bn$, since these conditions ensure the existence of
the $p$-adic $\log_p$ and $\exp_p$.

Denote
$$
S=\{\h=(h_n)_{n\in\bn}: \ h_n \ \textrm{satisfies} \
(\ref{eq1})\}.
$$

Let us rewrite (\ref{eq1}) as follows
\begin{equation}\label{eq2}
u_n=\left(\dsf{a_{n+1}u_{n+1}+b_{n+1}}{c_{n+1}+u_{n+1}}\right)
\left(\dsf{a_{n+2}u_{n+2}+b_{n+2}}{c_{n+2}+u_{n+2}}\right), \ \
n\in\bn,
\end{equation}
here we have denoted

\begin{equation}\label{d}
\left. \ba{ll}
 a_{k}=\exp_p(\a_{k}), \ \ \ b_{k}=\exp_p(\b_{k}), \\[3mm]
c_{k}=\exp_p(\g_{k}), \ \ \ \ u_{k}=\exp_p(h_{k}).\\[2mm]
\ea \right\}, \ \ k\in\bn.
\end{equation}

Let us introduce
\begin{equation}\label{f}
f_k(x)=\dsf{a_{k}x+b_{k}}{c_{k}+x}, \ \ k\in\bn.
\end{equation}

The functions $f_k$ have the the following properties.

\begin{lem}\label{2} For every $k\in\bn$ the following
relations hold
\begin{equation}\label{f1}
|f_k(x)-f_k(y)|_p\leq\frac{1}{p}\frac{|x-y|_p}{|c_k+x|_p|c_k+y|_p}
\end{equation}

\begin{equation}\label{f2}
|f_k(\exp_p(h))|_p=1 \ \ \textrm{for every} \ |h|_p\leq\frac{1}{p}
\end{equation}
\end{lem}

\begin{pf} From (\ref{f}) we immediately obtain
$$
|f_k(x)-f_k(y)|_p=\frac{|x-y|_p|a_kc_k-b_k|_p}{|c_k+x|_p|c_k+y|_p}.
$$
Keeping in mind (\ref{d}) from (\ref{exp}),(\ref{con}) one gets
$$
|a_kc_k-b_k|_p=|\exp_p(\a_k+\g_k-\b_k)-1|_p\leq\frac{1}{p}
$$
which implies (\ref{f1}). The strong triangle inequality with
(\ref{exp}) and $p\geq 3$ implies that $|\exp_p(h)+1|_p=1$, hence
one finds (\ref{f2}). \end{pf}

The main result of the paper is the following

\begin{thm}\label{3} Let (\ref{con}) be satisfied. Then
$|S|\leq 1$, here $|A|$ means the number of elements of a set $A$.
\end{thm}

\begin{pf} If $S=\emptyset$, then nothing to prove. So,
assume that $S\neq\emptyset$. To prove Theorem it is enough to
prove that any two elements of $S$ coincide with each other. To
show this it is enough to prove that for arbitrary $\varepsilon>0$
and every $\h=(h_n, n\in \bn), \z=(s_n,n\in\bn)\in S$ and
$n\in\bn$ the inequality $|h_n-s_n|_p<\varepsilon$ is valid.

Let $\h=(h_n,n\in \bn), \z=(s_n,n\in\bn)\in S$ and $\e>0$ be an
arbitrary number. Put $v_k=\exp_p(s_k)$. Let $n\in \bn$ be an
arbitrary number. Using (\ref{f}), (\ref{eq2}) and Lemma \ref{2}
we get
\begin{eqnarray}\label{ineq1}
|u_n-v_n|_p&= &
|f_{n+1}(u_{n+1})f_{n+2}(u_{n+2})-f_{n+1}(v_{n+1})f_{n+2}(v_{n+2})|_p\nonumber
\\
&\leq &
\max\bigg\{|f_{n+1}(u_{n+1})|_p|f_{n+2}(u_{n+2})-f_{n+2}(v_{n+2})|_p,
\nonumber \\
&&
|f_{n+2}(v_{n+2})|_p|f_{n+1}(u_{n+1})-f_{n+1}(v_{n+1})|_p\bigg\}
\nonumber \\
&\leq &\frac{1}{p}\max\{|u_{n+1}-v_{n+1}|_p,|u_{n+2}-v_{n+2}|_p\},
\end{eqnarray}
here we have used the following equalities
$$
|c_k+u_k|_p=1, \ \ |c_k+v_k|_p=1, \ \ \forall k\in\bn,
$$
which immediately follow from (\ref{exp}).

Now  take $n_0\in\bn$ such that $\dsf{1}{p^{n_0}}<\varepsilon$.
Then iterating (\ref{ineq1}) one gets
\[
|u_n-v_n|_p\leq\frac{1}{p^{n_0}}<\e
\]

The last inequality with (\ref{exp}), (\ref{d}) implies that
$$
|h_x-s_x|_p=|u_n-v_n|_p<\e
$$
This completes the proof. \end{pf}

{\bf Remark.} Note that the proved Theorem provides some
applications to the $p$-adic Bethe lattice models and models on
random tree-like graphs. This would be a theme of our further
investigations. We also think that the result can be applied for
the study of dynamics of random rational $p$-adic dynamical
systems. It should be also noted that certain type of recursive
equations in $p$-adic numbers were considered in \cite{M}, they
were related with transcendentality of the $p$-adic numbers. On
the other hand, the method presented here used functional approach
which is different from \cite{M}, where some number theoretical
method was presented. \\

Now consider several cases when $S$ is non empty.\\

{\it Case 1.} Let us assume that the following equality holds
\begin{equation*}
a_k+b_k=c_k+1
\end{equation*}
for every $k\in\bn$, which implies that the equation (\ref{eq2})
has a solution $u_n=1$ for all $n\in\bn$. According to Theorem
\ref{3} it is unique, i.e. $|S|=1$.\\

{\it Case 2.} Now suppose that $a_k=a$, $b_k=b$, $c_k=c$ for all
$k\in\bn$. We are going to show the existence of a solution
(\ref{eq2}). Let us search a solution in a form $u_n=u$, $\forall
n\in\bn$, here $u$ is unknown. Then (\ref{eq2}) reduces to the
equation  $f(u)=u$, where
$$
f(u)=\left(\dsf{a u+b}{c+u}\right)^2.
$$

Put $D=S(0,1)\cap B(1,p^{-1/(p-1)})$. Let us show that
$f(D)\subset D$. Indeed, let $x\in D$, then $|x|_p=1$,
$|x-1|_p<p^{-1/(p-1)}$. According to Lemma \ref{1} this means that
$x=\exp_p(y)$ for some $y\in B(0,p^{-1/(p-1)})$. Now by means of
Lemma \ref{2} we infer that $|f(x)|_p=1$. The strong triangle
inequality with
$$
|a-1|_p\leq\frac{1}{p}, \ \ |b-c|\leq \frac{1}{p}
$$
implies that
\begin{eqnarray}\label{ineq2}
|f(x)-1|_p&=&\bigg|\dsf{(a-1)x+b-c}{c+x}\bigg|_p\bigg|\dsf{(a-1)x+b+c}{c+x}\bigg|_p\leq
\frac{1}{p}.
\end{eqnarray}
Using the same argument along the proof of Lemma \ref{2} one gets
\begin{equation}\label{eq5}
|f(x)-f(y)|_p\leq\dsf{1}{p}|x-y|_p \ \ \textrm{for every} \ \
x,y\in D.
\end{equation}
Thus from the inequalities \eqref{ineq2}, \eqref{eq5} we obtain
that $f$ is a contraction of $D$, hence $f$ has a unique fixed
point $\zeta\in D$. Consequently, $S$ is non empty. Now Theorem
\ref{3} yields that $|S|=1$.

{\bf Acknowledgement}. The author thanks the FCT (Portugal) grant
SFRH/BPD/17419/2004, and also express gratitude to referees for
useful suggestions.

\end{document}